\documentclass[11pt]{amsart}
\usepackage{amssymb, latexsym}
\theoremstyle{plain}
\newtheorem{theorem}{Theorem}

\newtheorem*{1'}{Theorem 1-Bessel}
\newtheorem*{P2'}{Proposition 2-Bessel}
\newtheorem*{P3'}{Proposition 3-Bessel}
\newtheorem*{P4'}{Proposition 4-Bessel}
\newtheorem*{C1'}{Corollary 1-Bessel}

\newtheorem*{2'}{Theorem 2-Bessel}
\newtheorem*{3'}{Theorem 3-Bessel}

\theoremstyle{remark}

\newtheorem*{Remark 1}{Remark 1}
\newtheorem*{Remark 2}{Remark 2}
\newtheorem*{Remark 3}{Remark 3}
\newtheorem*{Remark 4}{Remark 4}

\numberwithin{equation}{section}

\begin{document}

\title [Inversions  in permutations with prescribed fixed points]
{The inversion statistic in derangements and in other permutations with a prescribed number of fixed points }

\author{Ross G. Pinsky}

\address{Department of Mathematics\\
Technion---Israel Institute of Technology\\
Haifa, 32000\\ Israel}
\email{ pinsky@technion.ac.il}

\urladdr{https://pinsky.net.technion.ac.il/}

\date{}

\begin{abstract}

We study how the inversion statistic is influenced by fixed points in a permutation.
For each $n\in\mathbb{N}$, and each $k\in\{0,1,\cdots, n\}$, let $P_n^{(k)}$ denote the uniform probability measure on the set of permutations in $S_n$ with exactly
$k$ fixed points. We obtain an exact formula for the expected number of inversions under the measure $P_n^{(k)}$ as well as for
$P_n^{(k)}(\sigma^{-1}_i<\sigma^{-1}_j)$, for $1\le i<j\le n$,  the $P_n^{(k)}$-probability that the number $i$ precedes the number $j$.
In particular,
 up  to a super-exponentially small correction as $n\to\infty$, the expected number of inversions in a random derangement $(k=0)$ is $\frac16n+\frac1{12}$ more than
 the value $\frac{n(n-1)}4$ that one obtains for a uniformly random
 general permutation in $S_n$.
On the other hand,  up  to a super-exponentially small correction,
for $k\ge2$, the expected number of inversions in  a random permutation with $k$ fixed points is $\frac{k-1}6n+\frac{k^2-k-1}{12}$ less than $\frac{n(n-1)}4$.
In  the borderline case, $k=1$, up to a
super-exponentially small correction, the expected number of inversions in a random permutation  with one fixed point is $\frac1{12}$ more than $\frac{n(n-1)}4$.
The proofs make strategic and perhaps novel use of the Chinese restaurant construction for a uniformly random permutation.
\end{abstract}

\maketitle
\section{Introduction and Statement of Results}\label{intro}
In this paper, we study how the inversion statistic is influenced by fixed points in a permutation. Let $P_n$ denote the uniform probability measure on the set $S_n$ of permutations of $[n]=\{1,\cdots, n\}$, and denote the expectation with respect to $P_n$ by $E_n$.
We write $\sigma\in S_n$ in one-line notation as $\sigma=\sigma_1\cdots\sigma_n$.
Recall that $\sigma$ is a \it derangement\rm\ if it has no fixed points; that is, if $\sigma_i\neq i$, for all $i\in[n]$.
Let $D_n\subset S_n$ denote the set of derangements in $S_n$.
As is well-known \cite{ABT, Pin14},
\begin{equation}\label{derang}
P_n(D_n)=\sum_{l=0}^n\frac{(-1)^l}{l!};\ \ \lim_{n\to\infty}P_n(D_n)=e^{-1}.
\end{equation}

Let $I_n(\sigma)$ denote the number of inversions in $\sigma$; that is,
$$
I_n(\sigma)=\sum_{1\le i<j\le n}1_{\{\sigma_j<\sigma_i\}}=\sum_{1\le i<j\le n}1_{\{\sigma^{-1}_j<\sigma^{-1}_i\}}.
$$
(Note that $\sigma^{-1}_i$ is the position of the number $i$ in the permutation $\sigma$.)
By symmetry, one has $E_nI_n=\frac{n(n-1)}4$.
It is well-known that the law of large numbers holds for $I_n$ under the uniform probability measure $P_n$; that is,
\begin{equation}\label{WLLN}
\lim_{n\to\infty} P_n\left((1-\epsilon)\frac{n^2}4\le I_n\le (1+\epsilon)\frac{n^2}4\right)=1,\ \text{for all}\ \epsilon>0.
\end{equation}

Denote the uniform probability measure on $D_n$ by $P_n^{(0)}$; that is,
$$
P_n^{(0)}(A)=\frac{P_n(A\cap D_n)}{P_n(D_n)},\  A\subset S_n,
$$
  and let $E_n^{(0)}$ denote the expectation with respect to $P_n^{(0)}$.
From \eqref{derang} and \eqref{WLLN} it follows that the law of large numbers also holds for $I_n$ under the measure $P_n^{(0)}$; that is, \eqref{WLLN} also holds with $P_n$ replaced by $P_n^{(0)}$.
From this and the fact that $\max_{\sigma\in S_n}I_n(\sigma)=O(n^2)$, if follows that  $E_n^{(0)}I_n\sim \frac{n^2}4$.  What about lower order terms in $E_n^{(0)}I_n$?
Intuitively, it is easy to see  that $E_n^{(0)}I_n$ should be larger than $E_nI_n=\frac{n(n-1)}4$. Indeed, by symmetry considerations, for any $i\in[n]$, one has $P_n^{(0)}(\sigma^{-1}_i=l)=\frac1{l-1},\ l\in[n]-\{i\}$.
Thus under $P_n^{(0)}$, for $1\le i<j\le n$, the random variable $\sigma^{-1}_i$ strictly stochastically dominates the random variable $\sigma^{-1}_j$; that is, 
$P_n^{(0)}(\sigma^{-1}_i\ge a)\ge P_n^{(0)}(\sigma^{-1}_j\ge a)$, for all $a\in \mathbb{R}$, and with strict inequality for at least one choice of $a$.
Of course this doesn't prove anything rigorously about $E_n^{(0)}I_n$ because  $\sigma^{-1}_i$ and $\sigma^{-1}_j$  are not independent.
We will calculate $E_n^{(0)}I_n$ and $P_n^{(0)}(\sigma^{-1}_i<\sigma^{-1}_j)$ explicitly, allowing one to   see how much greater $E_n^{(0)}I_n$ is than $E_nI_n$ and how much smaller
$P_n^{(0)}(\sigma^{-1}_i<\sigma^{-1}_j)$ is than   $P_n(\sigma^{-1}_i<\sigma^{-1}_j)$.

We now turn to permutations with a prescribed non-zero number of fixed points.
For $k\in \mathbb{N}$, let $D_{n;k}$ denote the set of permutations in $S_n$ with exactly $k$ fixed points; that is
$$
D_{n;k}=\{\sigma\in S_n: |\{i\in[n]: \sigma(i)=i\}|=k\}.
$$
For convenience, define $D_{n;0}=D_n$.
It is well-known \cite{ABT,Gon} that under the uniform measure $P_n$,
\begin{equation}\label{Poisson}
\lim_{n\to\infty}P_n(D_{n;k})=\frac{e^{-1}}{k!},\ k=0,1,\cdots;
\end{equation}
equivalently, under $P_n$, the random variable that counts the number of fixed points converges in distribution to the Poisson distribution with parameter one.
Let $P_n^{(k)}$ denote the uniform probability measure on $D_{n;k}$; that is,
$$
P_n^{(k)}(A)=\frac{P_n(A\cap D_{n;k})}{P_n(D_{n;k})},\ A\subset S_n,
$$
 and let $E_n^{(k)}$ denote the expectation with respect to $P_n^{(k)}$. 
From \eqref{WLLN} and \eqref{Poisson}, it follows that \eqref{WLLN} also holds with $P_n$ replaced by $P_n^{(k)}$. Thus, as with the case $k=0$, we have $E_n^{(k)}I_n\sim\frac{n^2}4$.
Which do we expect to be larger,  $E_n^{(k)}I_n$ or $E_nI_n=\frac{n(n-1)}4$? It is instructive to consider the extreme case in which $k=n$. Note that $P_n^{(n)}$ is the $\delta$-measure on the identity permutation; thus, $E_n^{(n)}I_n=0$.
It is not hard to  see intuitively   that the larger $k$ is, the smaller $E_n^{(k)}I_n$ should be.
Thus,  there should  be a threshold value of $k$ (perhaps depending on $n$), so that for $k$ less than the threshold,  $E_n^{(k)}I_n>E_nI_n$ and for $k$ above the threshold, $E_n^{(k)}I_n<E_nI_n$.
Since under the uniform measure, the expected number of fixed points is easily seen to be  equal to one for all $n$, the  candidate $k=1$ is an intuitive choice for the threshold, at least for sufficiently large $n$.
We will calculate $E_n^{(k)}I_n$
 explicitly.
In particular, we will see that the above noted threshold is indeed $k=1$, and we will see what happens at the threshold $k=1$.

The first theorem below concerns derangements and the second one treats permutations with a prescribed non-zero number of fixed points.

\begin{theorem}\label{thm1}
\noindent Let $n\ge3$. Then

\noindent i.
\begin{equation}\label{thm1result1}
\begin{aligned}
&E_n^{(0)}I_n=\frac{n(n-1)}4+\frac16n+\frac1{12}+\frac{\frac{(-1)^n}{n!}}{\sum_{l=0}^n\frac{(-1)^l}{l!}}\left(\frac{n-1}{12}\right)=\\
&\frac{n(n-1)}4+\frac16n+\frac1{12}+O\left(\frac1{(n-1)!}\right).
\end{aligned}
\end{equation}

\noindent ii.
\begin{equation}\label{thm1result2}
\begin{aligned}
&P_n^{(0)}(\sigma^{-1}_i<\sigma^{-1}_j)=\frac12+\frac{\big(1-2(j-i)\big)n+2(j-i)}{2n(n-1)(n-2)}+
\frac{\frac{(-1)^n}{n!}}{\sum_{l=0}^n\frac{(-1)^l}{l!}}\thinspace\frac{2(j-i)-n}{2n(n-2)}=\\
&\frac12+\frac{\big(1-2(j-i)\big)n+2(j-i)}{2n(n-1)(n-2)}+O\left(\frac1{(n+1)!}\right).
\end{aligned}
\end{equation}
\end{theorem}
\medskip

\newpage
\begin{theorem}\label{thm2}
Let $n\ge3$ and let  $k\in\{1,2,\cdots, n\}$. Then

\noindent i.
\begin{equation}\label{thm2result1}
\begin{aligned}
&E_n^{(k)}I_n=\frac{n(n-1)}4-\frac{k-1}6n-\frac{k^2-k-1}{12}+\frac{\frac{(-1)^{n-k}}{(n-k)!}}{\sum_{l=0}^{n-k}\frac{-1)^l}{l!}}\left(\frac{n-k-1}{12}\right)=\\
&\frac{n(n-1)}4-\frac{k-1}6n-\frac{k^2-k-1}{12}+O\left(\frac1{(n-k-1)!}\right).
\end{aligned}
\end{equation}
\noindent ii.

\begin{equation}\label{thm2result2}
\begin{aligned}
&P_n^{(k)}(\sigma^{-1}_i<\sigma^{-1}_j)=\\
&\frac12+\frac1{2n(n-1)(n-2)}\Big(\left(2(k-1)(j-i)+k^2-3k+1\right)n-2(k^2-k-1)(j-i)\Big)+\\
&O\left(\frac1{(n-k)!} \right).
\end{aligned}
\end{equation}
\end{theorem}

\bf\noindent Remark.\rm\ The precise formula for the term $O\left(\frac1{(n-k)!} \right)$ in \eqref{thm2result2} can be found at the end of the paper in \eqref{thm2result2full}.
\medskip

We now make some observations and comments regarding the above theorems. We begin with the expectation of $I_n$ in part (i) of Theorems \ref{thm1} and \ref{thm2}.
From part (i) of Theorem \ref{thm1}, one sees that up  to a super-exponentially small correction as $n\to\infty$, the expected number of inversions in the case of a random derangement is $\frac16n+\frac1{12}$ more than it is for a random general permutation.
On the other hand, from part (i) of Theorem \ref{thm2}, one sees that up  to a super-exponentially small correction,
for $k\ge2$, the expected number of inversions in  a random permutation with $k$ fixed points is $\frac{k-1}6n+\frac{k^2-k-1}{12}$ less than it is for a random general permutation.
In particular, in both the case of derangements ($k=0$) and in the case $k\ge2$, this difference grows linearly in $n$. Now consider the borderline case, $k=1$. As we noted above, for all $n$,
the expected number of fixed points in a random general permutation is one. Theorem \ref{thm2} shows that for all $n$, up to a
super-exponentially small correction, the expected number of inversions in a random permutation  with one fixed point is $\frac1{12}$ more than it is for a random general permutation.

As we noted above, the random variable that counts the number of fixed points in a random general permutation of size $n$  converges in distribution as $n\to\infty$ to the Poisson distribution with parameter one.
Let $X$ be a random variable with this distribution. Then $EX=1$ and $EX^2=2$. Not surprisingly, if we substitute $X$ for $k$ in the expression
$\frac{k-1}6n+\frac{k^2-k-1}{12}$ and take the expectation, we obtain $E(\frac{X-1}6n+\frac{X^2-X-1}{12})=0$.

Note that since Theorem \ref{thm2} is true for all $k\in\{1,\cdots, n\}$,
one can let
$k$ depend on $n$: $k=k_n$ with $\lim_{n\to\infty}k_n=\infty$.
It follows from  \eqref{thm2result1} that as long as $\lim_{n\to\infty}\frac{k_n}n<1$, then
 asymptotically as $n\to\infty$, up  to a super-exponentially small correction the expected number of inversions in  a random permutation with $k_n$ fixed points is $\frac{k_n-1}6n+\frac{k_n^2-k_n-1}{12}$ less than it is for a random general permutation.
In particular, if $k_n=o(n)$, then the correction is on the order $nk_n$, so $E_n^{(k_n)}I_n\sim E_nI_n\sim\frac{n^2}4$, but if $k_n\sim cn$, for $c\in(0,1)$, then
$E_n^{(k_n)}I_n\sim\frac{3-2c-c^2}{12}n^2$.

We now consider the probability of the event $\{\sigma^{-1}_i<\sigma^{-1}_j\}$ in part (ii) of Theorems \ref{thm1} and \ref{thm2}.
Consider first part (ii) of Theorem \ref{thm1} and let $i$ and $j$ depend on $n$: $i=i_n, j=j_n$ with $i_n<j_n$.
Then it follows from \eqref{thm1result2} that
$$
0\le \frac12-P_n^{(0)}(\sigma^{-1}_{i_n}<\sigma^{-1}_{j_n})=\theta(\frac{j_n-i_n}{n^2}).
$$
Consider now part (ii) of Theorem \ref{thm2}.
For  $i=i_n,j=j_n$ with $i_n<j_n$,  and $k=k_n\ge2$, it follows from \eqref{thm2result2} that   as long as $\lim_{n\to\infty}\frac{k_n}n<1$, then
$$
0\le P_n^{(0)}(\sigma^{-1}_{i_n}<\sigma^{-1}_{j_n})-\frac12=\theta\left(\frac{k_n(j_n-i_n)\wedge k_n^2)}{n^2}\right).
$$
However, for $i=i_n,j=j_n$ with $i_n<j_n$,  and $k=1$, it follows from \eqref{thm2result2} that
$$
P_n^{(0)}(\sigma^{-1}_{i_n}<\sigma^{-1}_{j_n})=\frac12+\frac1{2n(n-1)(n-2)}\big(-n+2(j_n-i_n)\big)+O(\frac1{(n-1)!}).
$$
Thus the sign of $P_n^{(0)}(\sigma^{-1}_{i_n}<\sigma^{-1}_{j_n})-\frac12$ depends on whether $j_n-i_n>\frac n2$ or $j_n-i_n<\frac n2$.
\medskip

The proofs of the two theorems make strategic and perhaps novel use of the Chinese restaurant construction for a uniformly random permutation.
We note that if one substitutes $k=0$ in \eqref{thm2result1} and \eqref{thm2result2}, these formulas reduce to \eqref{thm1result1} and \eqref{thm1result2} respectively.
We have presented the two results separately because derangements are an important class of permutations and because the formulas are considerably simpler for derangements. Furthermore, although the same general method is used to prove the two theorems, there are some technical differences in the proofs, and the calculations are much shorter for Theorem \ref{thm1}.
Theorem \ref{thm1} is proved in section \ref{pfthm1} and Theorem \ref{thm2} is proved in section \ref{pfthm2}.

We end this section with a description of the Chinese restaurant construction that will be used in the proofs.
This construction simultaneously yields a uniformly random permutation
$\Sigma_n$ in $S_n$, for all $n$ \cite{Pitman06,Pin14}. Furthermore, the construction is consistent  in the sense that if one writes out
the permutation $\Sigma_n$ as the product of its cycles  and deletes the number $n$ from the cycle in which it appears, then
the resulting random permutation of $S_{n-1}$ is equal to $\Sigma_{n-1}$.

The construction works as follows. Consider a restaurant with an unlimited number of circular tables, each of which
has an unlimited number of seats. Person number 1 sits at a table.  Now for $n\ge1$, suppose that   persons
number 1 through $n$ have already been seated. Then person number $n+1$ chooses a seat as follows.
For each $j\in[n]$, with probability $\frac1{n+1}$, person number $n+1$ chooses to sit to the left of person number $j$.
Also, with probability $\frac1{n+1}$, person number $n+1$ chooses to sit at an unoccupied table.
Now for each $n\in\mathbb{N}$, the random permutation $\Sigma_n\in S_n$ is defined by $\Sigma_n(i)=j$, if
after the first $n$ persons have taken seats,
person number $j$ is seated to the left of person number $i$.

\section{Proof of Theorem \ref{thm1}}\label{pfthm1}
For the duration of the proof, $n,i,j$ are fixed, with $1\le i<j\le n$.
We begin by calculating $P_n^{(0)}(\sigma^{-1}_i<\sigma^{-1}_j)$. This will prove
 part (ii) and will also be fundamental for the proof of part (i).
We implement the Chinese restaurant construction,  described at the end of section \ref{intro}, to build a uniformly random permutation in $S_n$.
However, we make one change.  From the construction, it is clear that the $n$  persons can enter in any order we like, as we still obtain a uniformly random permutation in $S_n$. We let the number $j$ be the last of the $n$ numbers to be used. That is, in the language of the construction, \it\ person $j$ chooses a seat  last,  after all the other $n-1$  persons with numbers in $[n]-\{j\}$ have already chosen their seats.\rm\ We denote by $\Sigma_n=\Sigma_n(1),\cdots, \Sigma_n(n)$   the uniformly random permutation in $S_n$ obtained via this construction.
Note that after $n-1$ stages of the construction, a uniformly random permutation of $[n]-\{j\}$ has been built.
We denote this permutation by
$\Sigma_{n-1}=\Sigma_{n-1}(1),\cdots, \Sigma_{n-1}(j-1),\Sigma_{n-1}(j+1),\cdots, \Sigma_{n-1}(n)$.
A fixed point for $\Sigma_{n-1}$ is a number $l\in[n]-\{j\}$ for which $\Sigma_{n-1}(l)=l$.
We note for later use that the distribution of the number of fixed points in
$\Sigma_{n-1}$ is the same as it is for a uniformly random permutation of $[n-1]$.
We use the generic $P$ to denote probabilities with respect to the above construction.

We now define several   events relative to the above construction. We note that we are
 are reserving the notation $D_l$ for the subset of derangements in $S_l, l\in\mathbb{N}$.
 Consequently, we  denote by
 $\mathcal{D}_{n-1}$  the event that $\Sigma_{n-1}$ has no fixed points and by
 $\mathcal{D}_n$ the event  that
 $\Sigma_n$ has no fixed points.
 Denote by $\mathcal{D}_{n-1,1}$ the event
 that $\Sigma_{n-1}$  has one fixed point.
Then from the construction, it follows that
\begin{equation}\label{Dnfromn-1}
\mathcal{D}_n=\left(\mathcal{D}_{n-1}\cap \mathcal{D}_n\right)\cup\left(\mathcal{D}_{n-1,1}\cap \mathcal{D}_n\right).
\end{equation}
Let
\begin{equation}\label{Cnij}
C_{n,i,j}=\{\Sigma^{-1}_n(i)<\Sigma^{-1}_n(j)\}.
\end{equation}
Note that $\Sigma^{-1}_n(i)$ denotes the position of $i$ in $\Sigma_n$.
By \eqref{Dnfromn-1},
\begin{equation}\label{DnCnij}
\mathcal{D}_n\cap C_{n,i,j}=\left(\mathcal{D}_{n-1}\cap \mathcal{D}_n\cap C_{n,i,j} \right)\cup\left(\mathcal{D}_{n-1,1}\cap \mathcal{D}_n\cap C_{n,i,j}\right).
\end{equation}
Thus,
\begin{equation}\label{translation}
P_n^{(0)}(\sigma^{-1}_i<\sigma^{-1}_j)=\frac{P(\mathcal{D}_{n-1}\cap \mathcal{D}_n\cap C_{n,i,j})}{P(\mathcal{D}_n)}+\frac{P(\mathcal{D}_{n-1,1}\cap \mathcal{D}_n\cap C_{n,i,j})}{P(\mathcal{D}_n)}.
\end{equation}

We now calculate the first term on the right hand side of \eqref{translation}.
From the construction, it follows that conditioned on $\mathcal{D}_{n-1}$, the random variable $\Sigma^{-1}_{n-1}(i)$ is uniformly distributed on $[n]-\{i,j\}$, and
conditioned on $\mathcal{D}_n$, the random variable $\Sigma^{-1}_n(j)$ is uniformly distributed on $[n]-\{j\}$. Furthermore, conditioned on $\mathcal{D}_{n-1}\cap \mathcal{D}_n$, the random variables
$\Sigma^{-1}_{n-1}(i)$ and $\Sigma^{-1}_n(j)$ are independent; thus the random vector $\left(\Sigma^{-1}_{n-1}(i),\Sigma^{-1}_n(j)\right)$ is uniformly distributed on
$\left([n]-\{i,j\}\right)\times\left( [n]-\{j\}\right)$.
Note from the construction that if $\left(\Sigma^{-1}_{n-1}(i),\Sigma^{-1}_n(j)\right)=(l_1,l_2)$ with $l_1\neq l_2$, then $\Sigma^{-1}_n(i)=l_1$, but if $l_1=l_2$, then $\Sigma^{-1}_n(i)=j$.
Thus, we conclude that
\begin{equation}\label{key1}
\begin{aligned}
&P\left(\left(\Sigma^{-1}_n(i),\Sigma^{-1}_n(j)\right)=(l_1,l_2)|\mathcal{D}_{n-1}\cap \mathcal{D}_n\right)=\\
&\begin{cases}\frac1{(n-2)(n-1)}; \ l_1\neq l_2, \ (l_1,l_2)\in\left([n]-\{i,j\}\right)\times\left( [n]-\{j\}\right);\\ \frac1{(n-2)(n-1)},\  l_1=j, l_2\in [n]-\{i,j\}.
\end{cases}.
\end{aligned}
\end{equation}

For each $m\in \mathbb{N}$, let $d_m$ denote the number of derangements in $S_m$. From \eqref{derang}, we have
\begin{equation}\label{countderang}
d_m=m!\left(\sum_{l=0}^m\frac{(-1)^l}{l!}\right).
\end{equation}
We write
\begin{equation}\label{PDn}
P(\mathcal{D}_m)=\frac{d_m}{m!}.
\end{equation}
From the construction, we have $P(\mathcal{D}_n|\mathcal{D}_{n-1})=\frac{n-1}n$. Using this with \eqref{PDn} gives
\begin{equation}\label{Dn-1Dn}
P(\mathcal{D}_{n-1}\cap \mathcal{D}_n)=(n-1)\frac{d_{n-1}}{n!}.
\end{equation}
From \eqref{Cnij}, \eqref{key1}, and \eqref{Dn-1Dn}, we conclude that
\begin{equation}\label{PDn-1DnC}
\begin{aligned}
&P(\mathcal{D}_{n-1}\cap \mathcal{D}_n\cap C_{n,i,j})=\\
&\frac{(n-1)d_{n-1}}{n!}\frac1{(n-2)(n-1)}\left(n-j+\sum_{l_1\in[j-1]-\{i\}}(n-1-l_1)+\sum_{l_1=j+1}^n(n-l_1)\right).
\end{aligned}
\end{equation}

 We now calculate the second term on the right hand side of \eqref{translation}.
 Let $d_{n-1,1}$ denote the number of permutations in $S_{n-1}$ with one fixed point.
 It follows easily that $d_{n-1,1}=(n-1)d_{n-2}$, where $d_{n-2}$ is as in \eqref{countderang}.
 So
 \begin{equation}\label{PDn-11}
 P(\mathcal{D}_{n-1,1})=\frac{d_{n-1,1}}{(n-1)!}=\frac{(n-1)d_{n-2}}{(n-1)!}=\frac{d_{n-2}}{(n-2)!}.
 \end{equation}
 Conditioned on $\mathcal{D}_{n-1,1}$, the event $\mathcal{D}_n$ occurs if and only if the number $j$, which enters at stage $n$, joins the lone singleton existing at stage $n-1$; the probability of this is of course $\frac1n$.
 Thus,
 \begin{equation}\label{Dn-11Dn}
 P(\mathcal{D}_n|\mathcal{D}_{n-1,1})=\frac1n.
 \end{equation}
 From  \eqref{PDn-11} and \eqref{Dn-11Dn}, we conclude that
 \begin{equation}\label{PDn-11Dn}
 P(\mathcal{D}_{n-1,1}\cap\mathcal{D}_n)=\frac{d_{n-2}}{n(n-2)!}.
 \end{equation}

 Now consider $P(C_{n,i,j}|\mathcal{D}_{n-1,1}\cap \mathcal{D}_n)$.
As we've already noted,  for the event $\mathcal{D}_{n-1,1}\cap \mathcal{D}_n$ to occur, the number $j$ must join the lone singleton existing at stage $n-1$. This singleton has equal probability of being
any number in $[n]-\{j\}$. In particular, with probability $\frac1{n-1}$, the number $j$ will join the number $i$. In this case $\Sigma^{-1}_n(i)=j$ and $\Sigma^{-1}_n(j)=i$, and thus the event
$\{\Sigma^{-1}_n(i)<\Sigma^{-1}_n(j)\}$ does not occur. Conditioned on  the singleton not being $i$,   it follows by symmetry that $(\Sigma^{-1}_n(i),\Sigma^{-1}_n(j))$  is distributed uniformly on
$\big([n]-\{i,j\}\big)\times\big([n]-\{i,j\}\big)$. Thus, in this case, the event $\{\Sigma^{-1}_n(i)<\Sigma^{-1}_n(j)\}$ occurs with probability $\frac12$.
From this, we conclude that
\begin{equation}\label{key2}
P(C_{n,i,j}|\mathcal{D}_{n-1,1}\cap \mathcal{D}_n)=\frac12\thinspace\frac{n-2}{n-1}.
\end{equation}
From \eqref{PDn-11Dn} and \eqref{key2}, we obtain
\begin{equation}\label{PDn-11DnC}
P(\mathcal{D}_{n-1,1}\cap \mathcal{D}_n\cap C_{n,i,j})=\frac12\thinspace\frac{n-2}{n-1}\frac{d_{n-2}}{n(n-2)!}=\frac12\frac{(n-2)d_{n-2}}{n!}.
\end{equation}

From \eqref{translation}, \eqref{PDn}, \eqref{PDn-1DnC} and \eqref{PDn-11DnC}, we obtain
\begin{equation}\label{probijfirst}
\begin{aligned}
&\frac{d_n}{n!}P_n^{(0)}(\sigma^{-1}_i<\sigma^{-1}_j)=\frac12\frac{(n-2)d_{n-2}}{n!}+\\
&\frac{(n-1)d_{n-1}}{n!}\frac1{(n-2)(n-1)}\left(n-j+\sum_{l_1\in[j-1]-\{i\}}(n-1-l_1)+\sum_{l_1=j+1}^n(n-l_1)\right).
\end{aligned}
\end{equation}
Using the fact that
$\sum_{l_1\in[j-1]-\{i\}}l_1=\frac{(j-1)j}2-i$ and that $\sum_{l_1=j+1}^nl_1=\frac{n(n+1)}2-\frac{j(j+1)}2$,
we can rewrite
\eqref{probijfirst} as
\begin{equation}\label{probijsecond}
\begin{aligned}
&\frac{d_n}{n!}P_n^{(0)}(\sigma^{-1}_i<\sigma^{-1}_j)=\frac12\frac{(n-2)d_{n-2}}{n!}+\\
&\frac{d_{n-1}}{(n-2)n!}\left((n-1)(j-2)+(n+1)(n-j)-\frac{(j-1)j}2+i-\frac{n(n+1)}2+\frac{j(j+1)}2\right).
\end{aligned}
\end{equation}
One has
$$
(n-1)(j-2)+(n+1)(n-j)-\frac{(j-1)j}2+i-\frac{n(n+1)}2+\frac{j(j+1)}2=\frac12n^2-\frac32n+i-j+2.
$$
Thus, from \eqref{probijsecond} we have
\begin{equation}\label{probijthird}
\begin{aligned}
&\frac{d_n}{n!}P_n^{(0)}(\sigma^{-1}_i<\sigma^{-1}_j)=\frac12\frac{(n-2)d_{n-2}}{n!}+
\frac{d_{n-1}}{(n-2)n!}\left( \frac12n^2-\frac32n+i-j+2 \right)=\\
&=\frac{d_{n-2}}{(n-2)!}\left(\frac12\frac{n-2}{n(n-1)}\right)+ \frac{d_{n-1}}{(n-1)!}\frac1{n(n-2)}\left( \frac12n^2-\frac32n+i-j+2 \right).
\end{aligned}
\end{equation}
From \eqref{countderang}, we have
\begin{equation}\label{dndn-1}
\frac{\frac{d_{n-1}}{(n-1)!}}{\frac{d_n}{n!}}=1-\frac{\frac{(-1)^n}{n!}}{\sum_{l=0}^n\frac{(-1)^l}{l!}},\ \ \frac{\frac{d_{n-2}}{(n-2)!}}{\frac{d_n}{n!}}=1-\frac{\frac{(-1)^{n-1}}{(n-1)!}+\frac{(-1)^n}{n!}}{\sum_{l=0}^n\frac{(-1)^l}{l!}}.
\end{equation}
Also, we have
$$
\frac1{n(n-2)}\left(\frac{n^2}2-\frac32n\right)=\frac12-\frac1{2(n-2)}.
$$
Using this with \eqref{dndn-1}, we can rewrite \eqref{probijthird} as
\begin{equation}\label{probijfourth}
\begin{aligned}
&P_n^{(0)}(\sigma^{-1}_i<\sigma^{-1}_j)=\left(1-\frac{\frac{(-1)^n}{n!}}{\sum_{l=0}^n\frac{(-1)^l}{l!}}\right)\left(\frac12-\frac1{2(n-2)}+\frac{2-(j-i)}{n(n-2)}\right)+\\
&\left(1-\frac{\frac{(-1)^{n-1}}{(n-1)!}+\frac{(-1)^n}{n!}}{\sum_{l=0}^n\frac{(-1)^l}{l!}}\right)\frac{n-2}{2n(n-1)}.
\end{aligned}
\end{equation}
After some algebra, which we leave to the reader, \eqref{probijfourth} can be written as given in \eqref{thm1result2}.
This proves part (ii) of the theorem.

We now turn to part (i) of the theorem.
We have
\begin{equation}\label{invernoninver}
E_n^{0}I_n=\frac{n(n-1)}2-\sum_{1\le i<j\le n}P_n^{(0)}(\sigma^{-1}_i<\sigma^{-1}_j).
\end{equation}
Using \eqref{thm1result2}, we calculate
$\sum_{1\le i<j\le n}P_n^{(0)}(\sigma^{-1}_i<\sigma^{-1}_j)$.
We note that
\begin{equation}\label{sumij}
\sum_{1\le i<j\le n}i=\frac16(n-1)n(n+1);\ \ \sum_{1\le i<j\le n}j=\frac13(n-1)n(n+1).
\end{equation}
We first sum over $1\le i<j\le n$ the terms on the right hand side of the first line of \eqref{thm1result2} that don't involve
$\frac{\frac{(-1)^n}{n!}}{\sum_{l=0}^n\frac{(-1)^l}{l!}}$. Using \eqref{sumij},  we have
\begin{equation}\label{firstsum}
\begin{aligned}
&\sum_{1\le i<j\le n}\left(\frac12+\frac{\big(1-2(j-i)\big)n+2(j-i)\Big)}{2n(n-1)(n-2)}\right)=\\
&\frac14n(n-1)-\frac n{2n(n-1)(n-2)}\frac12n(n-1)-\frac{2(n-1)}{2n(n-1)(n-2)}\frac16(n-1)n(n+1)=\\
& \frac14n(n-1)-\frac n6-\frac1{12}.
\end{aligned}
\end{equation}
Now we sum over $1\le i<j\le n$ the term multiplying
$\frac{\frac{(-1)^n}{n!}}{\sum_{l=0}^n\frac{(-1)^l}{l!}}$ on the right hand side of the first  line of \eqref{thm1result2}. Using \eqref{sumij}, we have
\begin{equation}\label{secondsum}
\sum_{1\le i<j\le n}\frac{2(j-i)-n}{2n(n-2)}=\frac1{n(n-2)}\frac16(n-1)n(n+1)-\frac 1{2(n-2)}\frac12n(n-1)=\frac{1-n}{12}.
\end{equation}
From \eqref{thm1result2}, \eqref{firstsum} and \eqref{secondsum}, we conclude that
\begin{equation}\label{sumPsigmaij}
\sum_{1\le i<j\le n}P_n^{(0)}(\sigma^{-1}_i<\sigma^{-1}_j)= \frac14n(n-1)-\frac n6-\frac1{12}+\frac{\frac{(-1)^n}{n!}}{\sum_{l=0}^n\frac{(-1)^l}{l!}}\left(\frac{1-n}{12}\right).
\end{equation}
Now \eqref{thm1result1} follows from \eqref{invernoninver} and \eqref{sumPsigmaij}. This completes the proof of part (i)
\hfill$\square$

\section{Proof of Theorem \ref{thm2}}\label{pfthm2}
For the duration of the proof, $n,i,j$ and $k$ are fixed,  with $1\le i<j\le n$ and $1\le k\le n$. We implement the Chinese restaurant construction in the same way that it was implemented in the proof of Theorem \ref{thm1}; that is, the number $j$ enters at the final stage.
 Recall from the first paragraph of the proof of Theorem \ref{thm1} that the uniformly random permutation in $S_n$ obtained via the construction is denoted by $\Sigma_n=\Sigma_n(1),\cdots, \Sigma_n(n)$, and the uniformly random permutation of $[n]-\{j\}$ obtained after the first $n-1$ stages of the construction is denoted by
$\Sigma_{n-1}=\Sigma_{n-1}(1),\cdots, \Sigma_{n-1}(j-1),\Sigma_{n-1}(j+1),\cdots, \Sigma_{n-1}(n)$.

We now define several   events relative to the above construction. We note that we are
 are reserving the notation $D_{n,l}$ for the subset of permutations  in $S_n$ with $l$ fixed points, $l\in[n]$.
 Consequently, we  denote by    $\mathcal{D}_{n,k}$ the event that $\Sigma_n$ has $k$ fixed points, and we denote by
 $\mathcal{D}_{n-1,l}$ the event that $\Sigma_{n-1}$ has $l$ fixed points, $l\in\{0,1,\cdots, n-1\}$.
From the construction, it follows that
\begin{equation}\label{Dnkfromn-1}
\mathcal{D}_{n,k}=\left(\mathcal{D}_{n-1,k-1}\cap\mathcal{D}_{n,k}\right)\cup\left(\mathcal{D}_{n-1,k}\cap\mathcal{D}_{n,k}\right)\cup\left(\mathcal{D}_{n-1,k+1}\cap\mathcal{D}_{n,k}\right).
\end{equation}
(If $k=n$, we understand $\mathcal{D}_{n-1,k}$ to be the empty set and if $k\in\{n-1,n\}$, we understand $\mathcal{D}_{n-1,k+1}$ to be the empty set.)

Let
$$
C_{n,i,j}=\{\Sigma^{-1}_n(i)<\Sigma^{-1}_n(j)\},
$$
as was defined in \eqref{Cnij}.
By \eqref{Dnkfromn-1}, we have
\begin{equation}\label{DnkCnij}
\begin{aligned}
&\mathcal{D}_{n,k}\cap C_{n,i,j}=\\
&\left(\mathcal{D}_{n-1,k-1}\cap\mathcal{D}_{n,k}\cap C_{n,i,j}\right)\cup\left(\mathcal{D}_{n-1,k}\cap\mathcal{D}_{n,k}\cap C_{n,i,j}\right)\cup\left(\mathcal{D}_{n-1,k+1}\cap\mathcal{D}_{n,k}\cap C_{n,i,j}\right).
\end{aligned}
\end{equation}
We will calculate the probability of each of the three events in the union on the right hand side of \eqref{DnkCnij}.

We begin with the  calculation of  $P(\mathcal{D}_{n-1,k-1}\cap\mathcal{D}_{n,k}\cap C_{n,i,j})$.
The event that  $i$ is a fixed point at stage $n-1$ is the event $\{\Sigma^{-1}_{n-1}(i)=i\}$. We  write
\begin{equation}\label{Dn-1k-1}
\begin{aligned}
&\mathcal{D}_{n-1,k-1}\cap\mathcal{D}_{n,k}=\\
&\big(\mathcal{D}_{n-1,k-1}\cap\mathcal{D}_{n,k}\cap\{\Sigma^{-1}_{n-1}(i)=i\}\big)\cup\big(\mathcal{D}_{n-1,k-1}\cap\mathcal{D}_{n,k}\cap\{\Sigma^{-1}_{n-1}(i)\neq i\}\big).
\end{aligned}
\end{equation}
Let $d_{m,l}$ denote the number of permutations in $S_m$ with $l$ fixed points.
Then $P(\mathcal{D}_{n-1,k-1})=\frac{d_{n-1,k-1}}{(n-1)!}$.
Given $\mathcal{D}_{n-1,k-1}$, the probability that $\Sigma^{-1}_{n-1}(i)=i$ is $\frac{k-1}{n-1}$.
Given $\mathcal{D}_{n-1,k-1}$, the event $\mathcal{D}_{n,k}$ will occur if and only if $j$ enters at stage $n$ as a fixed point; that is,
if and only if $\Sigma^{-1}_n(j)=j$. The probability of this is $\frac1n$.
On the event, $\mathcal{D}_{n-1,k-1}\cap\mathcal{D}_{n,k}\cap\{\Sigma^{-1}_{n-1}(i)=i\}$, one has $\Sigma^{-1}_n(i)=i$ and $\Sigma^{-1}_n(j)=j$, and thus the event $C_{n,i,j}$ occurs automatically.
From the above facts we obtain
\begin{equation}\label{first4intersect}
P\left(\mathcal{D}_{n-1,k-1}\cap\mathcal{D}_{n,k}\cap\{\Sigma^{-1}_{n-1}(i)=i\}\cap C_{n,i,j}\right)=\frac{d_{n-1,k-1}}{(n-1)!}\thinspace\frac{k-1}{n-1}\thinspace\frac1n.
\end{equation}
Given $\mathcal{D}_{n-1,k-1}$, the probability that $\Sigma^{-1}_{n-1}(i)\neq i$ is $\frac{n-k}{n-1}$.
Given the event $\mathcal{D}_{n-1,k-1}\cap\mathcal{D}_{n,k}\cap\{\Sigma^{-1}_{n-1}(i)\neq i\}$, $\Sigma^{-1}_n(i)$ is equally distributed over $[n]-\{i,j\}$. Thus, given this event, the event $C_{n,i,j}$ occurs
with probability $\frac{j-2}{n-2}$. From these facts  we obtain
\begin{equation}\label{second4intersect}
P\left(\mathcal{D}_{n-1,k-1}\cap\mathcal{D}_{n,k}\cap\{\Sigma^{-1}_{n-1}(i)\neq i\}\cap C_{n,i,j}\right)=\frac{d_{n-1,k-1}}{(n-1)!}\thinspace\frac{n-k}{n-1}\thinspace\frac1n\thinspace\frac{j-2}{n-2}.
\end{equation}
From \eqref{Dn-1k-1}--\eqref{second4intersect}, we conclude that
\begin{equation}\label{firstterm}
P(\mathcal{D}_{n-1,k-1}\cap\mathcal{D}_{n,k}\cap C_{n,i,j})=
\frac{d_{n-1,k-1}}{(n-1)!}\left(\frac{k-1}{n(n-1)}+\frac{(n-k)(j-2)}{n(n-1)(n-2)}\right).
\end{equation}

We now turn to the calculation of $P(\mathcal{D}_{n-1,k}\cap\mathcal{D}_{n,k}\cap C_{n,i,j})$.
Similar to what we did in \eqref{Dn-1k-1}, we write
\begin{equation}\label{similar}
\begin{aligned}
&\mathcal{D}_{n-1,k}\cap\mathcal{D}_{n,k}=\\
&\left(\mathcal{D}_{n-1,k}\cap\mathcal{D}_{n,k}\cap\{\Sigma^{-1}_{n-1}(i)=i\}\right)
\cup\left( \mathcal{D}_{n-1,k}\cap\mathcal{D}_{n,k}\cap\{\Sigma^{-1}_{n-1}(i)\neq i\}\right).
\end{aligned}
\end{equation}
We have $P(\mathcal{D}_{n-1,k})=\frac{d_{n-1,k}}{(n-1)!}$.
Given $\mathcal{D}_{n-1,k}$, the probability that $\Sigma^{-1}_{n-1}(i)=i$ is $\frac k{n-1}$.
Given $\mathcal{D}_{n-1,k}$, the event $\mathcal{D}_{n,k}$ will occur if and only if person $j$ entering  at stage $n$
sits at a table that already has at least two people. (If person $j$ starts a new table, then $\Sigma_n$ will have $k+1$ fixed points, and if person $j$ sits at a table that has one person, then $\Sigma_n$ will have $k-1$ fixed points.)
 The probability of this is $\frac{n-k-1}n$.
Given the event $\mathcal{D}_{n-1,k}\cap\mathcal{D}_{n,k}\cap\{\Sigma^{-1}_{n-1}(i)=i\}$, $\Sigma^{-1}_n(j)$ is uniformly distributed over $[n]-\{i,j\}$; thus
given this event, the probability that $C_{n,i,j}$ occurs is $\frac{n-i-1}{n-2}$.
From these facts, we obtain
\begin{equation}\label{Dnkfirst}
P\left(\mathcal{D}_{n-1,k}\cap\mathcal{D}_{n,k}\cap\{\Sigma^{-1}_{n-1}(i)=i\}\cap C_{n,i,j}\right)=\frac{d_{n-1,k}}{(n-1)!}\thinspace \frac k{n-1}\thinspace\frac{n-k-1}n\thinspace\frac{n-i-1}{n-2}.
\end{equation}

We now calculate $P\left(\mathcal{D}_{n-1,k}\cap\mathcal{D}_{n,k}\cap\{\Sigma^{-1}_{n-1}(i)\neq i\}\cap C_{n,i,j}\right)$.
Given $\mathcal{D}_{n-1,k}$, the probability that $\Sigma^{-1}_{n-1}(i)\neq i$ is $\frac {n-k-1}{n-1}$.
As noted above,
given $\mathcal{D}_{n-1,k}$, the event $\mathcal{D}_{n,k}$ will occur if and only if person $j$ entering  at stage $n$
sits at a table that already has at least two people, and
the probability of this is $\frac{n-k-1}n$. We now need to calculate the probability of $C_{n,i,j}$, given
$\mathcal{D}_{n-1,k}\cap\mathcal{D}_{n,k}\cap\{\Sigma^{-1}_{n-1}(i)\neq i\}$.
We need to be somewhat careful here. Conditioned on the above event, $\Sigma^{-1}_{n-1}(i)$ is uniformly distributed over
$[n]-\{i,j\}$. Conditioned on the above event, $\Sigma^{-1}_n(j)$ is uniformly distributed over the $n-k-1$ numbers that are not fixed points for $\Sigma_{n-1}$. These $n-k-1$ numbers include $i$.
The other $n-k-2$ such numbers are uniformly distributed over  all the $(n-k-2)$-tuples in   $[n]-\{i,j\}$.
From this, it follows that
\begin{equation}\label{lots}
\begin{aligned}
&P(\Sigma^{-1}_{n-1}(i)=l_1,\Sigma^{-1}_n(j)=i)=\frac1{n-2}\thinspace\frac1{n-k-1};\ l_1\in[n]-\{i,j\};\\
&P(\Sigma^{-1}_{n-1}(i)=l_1,\Sigma^{-1}_n(j)=l_1)=\frac1{n-2}\thinspace\frac1{n-k-1},\ l_1\in[n]-\{i,j\};\\
&P(\Sigma^{-1}_{n-1}(i)=l_1,\Sigma^{-1}_n(j)=l_2)=\frac1{n-2}\thinspace\frac{n-k-3}{n-k-1}\thinspace\frac1{n-3},\ l_1\neq l_2,\ l_1,l_2\in[n]-\{i,j\}.
\end{aligned}
\end{equation}
If $\{\Sigma^{-1}_{n-1}(i)=l_1,\Sigma^{-1}_n(j)=i\}$ occurs with $l_1\in[n]-\{i,j\}$, then also $\Sigma^{-1}_n(i)=l_1$, and if
$\{\Sigma^{-1}_{n-1}(i)=l_1,\Sigma^{-1}_n(j)=l_2)\}$ occurs with $l_1\neq l_2,\ l_1,l_2\in[n]-\{i,j\}$, then also $\Sigma^{-1}_n(i)=l_1$.
However, if $\{\Sigma^{-1}_{n-1}(i)=l_1,\Sigma^{-1}_n(j)=l_1\}$ occurs with $l_1\in[n]-\{i,j\}$, then from the construction we have
$\Sigma^{-1}_n(i)=j$. Using this with \eqref{lots}, it follows that conditioned on  $\mathcal{D}_{n-1,k}\cap\mathcal{D}_{n,k}\cap\{\Sigma^{-1}_{n-1}(i)\neq i\}$,
the probability of $C_{n,i,j}$ is $\frac{i-1}{n-2}\thinspace\frac1{n-k-1}+\frac{n-j}{n-2}\thinspace\frac1{n-k-1}+\frac12\frac{n-k-3}{n-k-1}$, which we write as
$\frac1{(n-2)(n-k-1)}\left(i-1+n-j+\frac12(n-2)(n-k-3)\right)$.
Thus, we conclude that
\begin{equation}\label{Dnksecond}
\begin{aligned}
&P\left(\mathcal{D}_{n-1,k}\cap\mathcal{D}_{n,k}\cap\{\Sigma^{-1}_{n-1}(i)\neq i\}\cap C_{n,i,j}\right)=\\
&\frac{d_{n-1,k}}{(n-1)!}\thinspace\frac{n-k-1}{n-1}\thinspace\frac{n-k-1}n
\frac1{(n-2)(n-k-1)}\left(i-1+n-j+\frac12(n-2)(n-k-3)\right).
\end{aligned}
\end{equation}
From \eqref{similar}, \eqref{Dnkfirst} and \eqref{Dnksecond} we obtain
\begin{equation}\label{secondterm}
\begin{aligned}
&P(\mathcal{D}_{n-1,k}\cap\mathcal{D}_{n,k}\cap C_{n,i,j})=\frac{d_{n-1,k}}{(n-1)!}\times\\
&\Big(\frac{k(n-k-1)(n-i-1)}{n(n-1)(n-2)}+
\frac{(n-k-1)}{n(n-1)(n-2)}\left(i-1+n-j+\frac12(n-2)(n-3-k)\right)\Big).
\end{aligned}
\end{equation}

We now turn to the calculation of $P\left(\mathcal{D}_{n-1,k+1}\cap\mathcal{D}_{n,k}\cap C_{n,i,j}\right)$.
Similar to before, we write
\begin{equation}\label{similaragain}
\begin{aligned}
&\mathcal{D}_{n-1,k+1}\cap\mathcal{D}_{n,k}=\\
&\left(\mathcal{D}_{n-1,k+1}\cap\mathcal{D}_{n,k}\cap\{\Sigma^{-1}_{n-1}(i)=i\}\right)
\cup\left( \mathcal{D}_{n-1,k+1}\cap\mathcal{D}_{n,k}\cap\{\Sigma^{-1}_{n-1}(i)\neq i\}\right).
\end{aligned}
\end{equation}
We have $P(\mathcal{D}_{n-1,k+1})=\frac{d_{n-1,k+1}}{(n-1)!}$.
Given $\mathcal{D}_{n-1,k+1}$, the probability that $\Sigma^{-1}_{n-1}(i)=i$ is $\frac{k+1}{n-1}$.
Given $\mathcal{D}_{n-1,k+1}$, the event $\mathcal{D}_{n,k}$ will occur if and only if person $j$ entering at stage $n$ sits at a table that has
exactly one person already.
(This reduces the number of fixed points from $k+1$ at stage $n-1$ to $k$ at stage $n$.) The probability of this is $\frac{k+1}n$.
We now need to calculate the probability of $C_{n,i,j}$, given
$\mathcal{D}_{n-1,k+1}\cap\mathcal{D}_{n,k}\cap\{\Sigma^{-1}_{n-1}(i)=i\}$. Given this event, person $j$ has probability $\frac1{k+1}$ of joining $i$ in which case
the event $C_{n,i,j}$ does not occur. With probability $\frac k{k+1}$, person $j$ will join  a table with a person other than $i$, and this person is uniformly distributed over
$[n]-\{i,j\}$. Thus given the above event, the probability of $C_{n,i,j}$ is $\frac k{k+1}\frac{n-1-i}{n-2}$.
Thus, we conclude that
\begin{equation}\label{Dnk+1first}
\begin{aligned}
&P\left(\mathcal{D}_{n-1,k+1}\cap\mathcal{D}_{n,k}\cap\{\Sigma^{-1}_{n-1}(i)= i\}\cap C_{n,i,j}\right)=\\
&\frac{d_{n-1,k+1}}{(n-1)!}\thinspace\frac{k+1}{n-1}\thinspace\frac{k+1}n\thinspace\frac k{k+1}\thinspace\frac{n-i-1}{n-2}.
\end{aligned}
\end{equation}

We now calculate $P\left(\mathcal{D}_{n-1,k+1}\cap\mathcal{D}_{n,k}\cap\{\Sigma^{-1}_{n-1}(i)\neq i\}\cap C_{n,i,j}\right)$.
Given $\mathcal{D}_{n-1,k+1}$, the probability that $\Sigma^{-1}_{n-1}(i)\neq i$ is $\frac{n-2-k}{n-1}$.
As noted above, given $\mathcal{D}_{n-1,k+1}$, the event $\mathcal{D}_{n,k}$ will occur if and only if person $j$ entering at stage $n$ sits at a table that has
exactly one person already, and the probability of this is $\frac{k+1}n$. We now need to calculate the probability of $C_{n,i,j}$, given
$\mathcal{D}_{n-1,k+1}\cap\mathcal{D}_{n,k}\cap\{\Sigma^{-1}_{n-1}(i)\neq i\}$. Given this event, $\Sigma^{-1}_n(i)$ and $\Sigma^{-1}_n(j)$ are independent and uniformly distributed over $[n]-\{i,j\}$; thus the probability of $C_{n,i,j}$ is $\frac12$.
Thus, we conclude that
\begin{equation}\label{Dnk+1second}
\begin{aligned}
&P\left(\mathcal{D}_{n-1,k+1}\cap\mathcal{D}_{n,k}\cap\{\Sigma^{-1}_{n-1}(i)\neq i\}\cap C_{n,i,j}\right)=\\
&\frac{d_{n-1,k+1}}{(n-1)!}\thinspace\frac{n-2-k}{n-1}\thinspace\frac{k+1}n\thinspace\frac12.
\end{aligned}
\end{equation}
From \eqref{similaragain}--\eqref{Dnk+1second}, we obtain
\begin{equation}\label{thirdterm}
\begin{aligned}
P(\mathcal{D}_{n-1,k+1}\cap\mathcal{D}_{n,k}\cap C_{n,i,j})=\frac{d_{n-1,k+1}}{(n-1)!}\left
(\frac{(k+1)k(n-i-1)}{n(n-1)(n-2)}+\frac{(n-2-k)(k+1)}{2n(n-1)}\right).
\end{aligned}
\end{equation}

From \eqref{DnkCnij}, \eqref{firstterm}, \eqref{secondterm} and \eqref{thirdterm},
it follows that $P_n(D_{n,k}\cap \{\sigma^{-1}_n(i)<\sigma^{-1}_n(j)\})=P(\mathcal{D}_{n,k}\cap C_{n,i,j})$ is equal to the sum of the terms on the right hand sides of
\eqref{firstterm}, \eqref{secondterm} and \eqref{thirdterm}.
Of course, $P_n^{(k)}(\sigma^{-1}_n(i)<\sigma^{-1}_n(j))=\frac1{P_n(D_{n,k})}P_n(D_{n,k}\cap \{\sigma^{-1}_n(i)<\sigma^{-1}_n(j)\})$.
Also
\begin{equation}\label{kinvernoninver}
E_n^{(k)}I_n=\frac{n(n-1)}2-\sum_{1\le i<j\le n}P_n^{(k)}(\sigma^{-1}_i<\sigma^{-1}_j).
\end{equation}
Thus to calculate $E_n^{(k)}I_n$, we will calculate
$$
\sum_{1\le i<j\le k}P_n(D_{n,k}\cap\{\sigma^{-1}_n(i)<\sigma^{-1}_n(j)\}),
$$
 which is the sum of
\eqref{firstterm}, \eqref{secondterm} and \eqref{thirdterm} over the pairs $i,j$ satisfying $1\le i<j\le n$.

Using \eqref{sumij}, we have for the sum  of \eqref{firstterm} over the pairs $i,j$ satisfying $1\le i<j\le n$,
\begin{equation}\label{firsttermsum}
\begin{aligned}
&\sum_{1\le i<j\le n}P(\mathcal{D}_{n-1,k-1}\cap\mathcal{D}_{n,k}\cap C_{n,i,j})=\\
&\sum_{1\le i<j\le n}\frac{d_{n-1,k-1}}{(n-1)!}\left(\frac{k-1}{n(n-1)}+\frac{(n-k)(j-2)}{n(n-1)(n-2)}\right)=\frac{d_{n-1,k-1}}{(n-1)!}\times\\
&\left(\frac{k-1}{n(n-1)}\frac{n(n-1)}2-\frac{2(n-k)}{n(n-1)(n-2)}\frac{n(n-1)}2+\frac{n-k}{n(n-1)(n-2)}\frac{(n-1)n(n+1)}3\right)=\\
&\frac{d_{n-1,k-1}}{(n-1)!}\left(\frac{k-1}2+\frac{n-k}3\right).
\end{aligned}
\end{equation}

We now consider the sum  of \eqref{secondterm} over the pairs $i,j$ satisfying $1\le i<j\le n$.
Using \eqref{sumij}, the sum of the first term on the second line of \eqref{secondterm} is
\begin{equation}\label{secondtermsum1}
\begin{aligned}
&\sum_{1\le i<j\le n}\frac{k(n-k-1)(n-i-1)}{n(n-1)(n-2)}=\\
&\frac{k(n-k-1)(n-1)}{n(n-1)(n-2)}\frac{n(n-1)}2-\frac{k(n-k-1)}{n(n-1)(n-2)}\frac{(n-1)n(n+1)}6=\\
&k\frac{n-k-1}{n-2}\left(\frac{n-1}2-\frac{n+1}6\right).
\end{aligned}
\end{equation}
Using \eqref{sumij}, the sum of the second term on the second line of \eqref{secondterm} is
\begin{equation}\label{secondtermsum2}
\begin{aligned}
&\sum_{1\le i<j\le n}\frac{(n-k-1)}{n(n-1)(n-2)}\left(i-1+n-j+\frac12(n-2)(n-3-k)\right)=\\
&\frac{n-k-1}{n(n-1)(n-2)}\left(n-1+\frac12(n-2)(n-3-k)\right)\frac{n(n-1)}2-\\
&\frac{n-k-1}{n(n-1)(n-2)}\frac{(n-1)n(n+1)}6=\frac{n-k-1}{n-2}\left(\frac{n-1}2+\frac{(n-2)(n-3-k)}4-\frac{n+1}6\right).
\end{aligned}
\end{equation}
Thus, from \eqref{secondterm}, \eqref{secondtermsum1} and \eqref{secondtermsum2}, we obtain
\begin{equation}\label{secondtermsum}
\begin{aligned}
&\sum_{1\le i<j\le n}P(\mathcal{D}_{n-1,k}\cap\mathcal{D}_{n,k}\cap C_{n,i,j})=\frac{d_{n-1,k}}{(n-1)!}\frac{n-k-1}{n-2}\times\\
&\left(\frac{k(n-1)}2-\frac{k(n+1)}6+\frac{n-1}2+\frac{(n-2)(n-3-k)}4-\frac{n+1}6\right)=\\
&\frac{d_{n-1,k}}{(n-1)!}\frac{n-k-1}{12(n-2)}\left(3n^2-11n+(n-2)k+10\right).
\end{aligned}
\end{equation}

Now we consider the sum of \eqref{thirdterm} over the pairs $i,j$ satisfying $1\le i<j\le n$.
Using \eqref{sumij}, we have
\begin{equation}\label{thirdtermsum}
\begin{aligned}
&\sum_{1\le i<j\le n}P(\mathcal{D}_{n-1,k+1}\cap\mathcal{D}_{n,k}\cap C_{n,i,j})=\\
&\sum_{1\le i<j\le n}\frac{d_{n-1,k+1}}{(n-1)!}\left
(\frac{(k+1)k(n-i-1)}{n(n-1)(n-2)}+\frac{(n-2-k)(k+1)}{2n(n-1)}\right)=\\
&\frac{d_{n-1,k+1}}{(n-1)!}
\Bigg(\left(\frac{(k+1)k(n-1)}{n(n-1)(n-2)}+\frac{(n-2-k)(k+1)}{2n(n-1)}\right)\frac{n(n-1)}2-\\
&\frac{(k+1)k}{n(n-1)(n-2)}\frac{(n-1)n(n+1)}6\Bigg)=\frac{d_{n-1,k+1}}{(n-1)!}\left(\frac{k+1}4n+\frac1{12}(k+1)(k-6)\right).
\end{aligned}
\end{equation}

We conclude that  $\sum_{1\le i<j\le n}P_n(D_n\cap\{\sigma^{-1}_n(i)<\sigma^{-1}_n(j)\})$,  is the sum of the right hand sides of \eqref{firsttermsum}, \eqref{secondtermsum} and \eqref{thirdtermsum}.
Consequently,  $\sum_{1\le i<j\le n}P^{(k)}_n(\sigma^{-1}_n(i)<\sigma^{-1}_n(j))$ is the above noted sum divided by $P_n(D_{n,k})$.
It is easy to see that  $d_{m,r}=\binom mrd_{m-r}$, for $0\le r\le m$. Thus,
\begin{equation}\label{dmr}
\frac{d_{m,r}}{m!}=\frac1{r!}\frac{d_{m-r}}{(m-r)!}=\frac1{r!}\sum_{l=0}^{m-r}\frac{(-1)^l}{l!},
\end{equation}
where the last equality follows from \eqref{countderang}.
For convenience in notation, let
\begin{equation}\label{Em-1}
\mathcal{E}_m(-1)=\sum_{l=0}^m\frac{(-1)^l}{l!}.
\end{equation}
Thus, from \eqref{dmr},
$$
P_n(D_{n,k})=\frac{d_{n,k}}{n!}=\frac1{k!}\sum_{l=0}^{n-k}\frac{(-1)^l}{l!}=\frac1{k!}\mathcal{E}_{n-k}(-1).
$$
Using this with \eqref{dmr}, we have
\begin{equation}\label{dnkquotients}
\begin{aligned}
&\frac{\frac{d_{n-1,k-1}}{(n-1)!}}{P_n(D_{n,k})}=k;\\
&\frac{\frac{d_{n-1,k}}{(n-1)!}}{P_n(D_{n,k})}=
1-\frac{\frac{(-1)^{n-k}}{(n-k)!}}{\mathcal{E}_{n-k}(-1)};\\
&\frac{\frac{d_{n-1,k+1}}{(n-1)!}}{P_n(D_{n,k})}=
\frac1{k+1}\left(1-\frac{\frac{(-1)^{n-k-1}}{(n-k-1)!}+\frac{(-1)^{n-k}}{(n-k)!}}{\mathcal{E}_{n-k}(-1)}   \right).
\end{aligned}
\end{equation}

From the facts noted in the previous paragraph along with \eqref{dnkquotients}, \eqref{firsttermsum}, \eqref{secondtermsum} and \eqref{thirdtermsum}, we obtain
\begin{equation*}
\begin{aligned}
&\sum_{1\le i<j\le n}P_n^{(k)}(\sigma^{-1}_n(i)<\sigma^{-1}_n(j))=k\left(\frac{k-1}2+\frac{n-k}3\right)+\\
&\left(1-\frac{\frac{(-1)^{n-k}}{(n-k)!}}{\mathcal{E}_{n-k}(-1)}\right)
\left(\frac{n-k-1}{12(n-2)}\left(3n^2-11n+(n-2)k+10\right)\right)+\\
&\frac1{k+1}\left(1-\frac{\frac{(-1)^{n-k-1}}{(n-k-1)!}+\frac{(-1)^{n-k}}{(n-k)!}}{\mathcal{E}_{n-k}(-1)}   \right)
\left(\frac{k+1}4n+\frac1{12}(k+1)(k-6)\right).
\end{aligned}
\end{equation*}
From this and \eqref{kinvernoninver}, we obtain
\begin{equation}\label{EI}
\begin{aligned}
&E_n^{(k)}I_n=\frac{n(n-1)}2-k\left(\frac{k-1}2+\frac{n-k}3\right)-\\
&\left(1-\frac{\frac{(-1)^{n-k}}{(n-k)!}}{\mathcal{E}_{n-k}(-1)}\right)
\left(\frac{n-k-1}{12(n-2)}\left(3n^2-11n+(n-2)k+10\right)\right)-\\
&\left(1-\frac{\frac{(-1)^{n-k-1}}{(n-k-1)!}+\frac{(-1)^{n-k}}{(n-k)!}}{\mathcal{E}_{n-k}(-1)}   \right)
\left(\frac n4+\frac1{12}(k-6)\right).
\end{aligned}
\end{equation}
Now
\begin{equation}\label{simplification}
\frac{n-k-1}{12(n-2)}(3n^2-11n+(n-2)k+10)=\frac14n^2-\frac{4+k}6n-\frac1{12}(k+1)(k-5).
\end{equation}
Consequently, the summands on the right hand side of \eqref{EI} that do not involve $\mathcal{E}_n(-1)$ are given by
\begin{equation}\label{notinvolve}
\begin{aligned}
&\frac{n(n-1)}2-k\left(\frac{k-1}2+\frac{n-k}3\right)
-\frac14n^2+\frac{4+k}6n+\frac1{12}(k+1)(k-5)-\frac n4-\frac1{12}(k-6)=\\
&\frac{n(n-1)}4-\frac{k-1}6n-\frac{k^2-k-1}{12}.
\end{aligned}
\end{equation}
Writing $\frac{(-1)^{n-k-1}}{(n-k-1)!}=-(n-k)\frac{(-1)^{n-k}}{(n-k)!}$ and using \eqref{simplification}, we can write the summands on the right hand side of \eqref{EI} that involve
$\mathcal{E}_n(-1)$ as
\begin{equation}\label{involve}
\begin{aligned}
&\frac{\frac{(-1)^{n-k}}{(n-k)!}}{\mathcal{E}_{n-k}(-1)}\left(\frac14n^2-\frac{4+k}6n-\frac1{12}(k+1)(k-5)+\big(-(n-k)+1\big)\big(\frac n4+\frac1{12}(k-6)\big)\right)=\\
&\frac{\frac{(-1)^{n-k}}{(n-k)!}}{\mathcal{E}_{n-k}(-1)}\left(\frac{n-k-1}{12}\right).
\end{aligned}
\end{equation}
Now \eqref{thm2result1} follows from \eqref{EI}, \eqref{notinvolve}, \eqref{involve} and \eqref{Em-1}.

We now prove \eqref{thm2result2}. From the first two sentences in the paragraph following \eqref{thirdterm}, along with \eqref{dnkquotients}, it follows that
\begin{equation}\label{Pijk1}
\begin{aligned}
&P_n^{(k)}(\sigma^{-1}_n(i)<\sigma^{-1}_n(j))=k\left(\frac{k-1}{n(n-1)}+\frac{(n-k)(j-2)}{n(n-1)(n-2)}\right)+\\
&\left(1-\frac{\frac{(-1)^{n-k}}{(n-k)!}}{\mathcal{E}_{n-k}(-1)}\right)\times\\
&\left(\frac{k(n-k-1)(n-i-1)}{n(n-1)(n-2)}+
\frac{(n-k-1)}{n(n-1)(n-2)}\bigg(i-1+n-j+\frac12(n-2)(n-3-k)\bigg)\right)+\\
&\frac1{k+1}\left(1-\frac{\frac{(-1)^{n-k-1}}{(n-k-1)!}+\frac{(-1)^{n-k}}{(n-k)!}}{\mathcal{E}_{n-k}(-1)}   \right)\times\\
&\left(\frac{(k+1)k(n-i-1)}{n(n-1)(n-2)}+\frac{(n-2-k)(k+1)}{2n(n-1)}\right).
\end{aligned}
\end{equation}
With a lot of algebra, one can show that the summands on the right hand side of \eqref{Pijk1} that do not involve $\mathcal{E}_n(-1)$ satisfy
\begin{equation}\label{Pijknotinvolve}
\begin{aligned}
&k\left(\frac{k-1}{n(n-1)}+\frac{(n-k)(j-2)}{n(n-1)(n-2)}\right)+\\
&\frac{k(n-k-1)(n-i-1)}{n(n-1)(n-2)}+  \frac{(n-k-1)}{n(n-1)(n-2)}\bigg(i-1+n-j+\frac12(n-2)(n-3-k)\bigg)+\\
&\frac{k(n-i-1)}{n(n-1)(n-2)}+\frac{(n-2-k)}{2n(n-1)}=\\
&\frac12+\frac1{2n(n-1)(n-2)}\left(\big(2(k-1)(j-i)+(k^2-3k+1)n-2(k^2-k-1)(j-i)\right).
\end{aligned}
\end{equation}
Now \eqref{thm2result2} follows from \eqref{Pijk1}, \eqref{Pijknotinvolve} and \eqref{Em-1}
\hfill $\square$

With considerably more algebra, one can show that the summands  on the  right hand side of \eqref{Pijk1} that involve $\mathcal{E}_n(-1)$ can be written as
\begin{equation}\label{canshow}
\begin{aligned}
&\frac{\frac{(-1)^{n-k}}{(n-k)!}}{\sum_{l=0}^{n-k}\frac{-1)^l}{l!}}\thinspace\frac1{2n(n-1)(n-2)}\times\\
&\Big(kn^3-(5k+1)n^2+(-k^3(3-2i)k^2+9k+2(j-i)+1)n+\\
&2ik^3+2(i-2)k^2-2(j-i+2)k-2(j-i)\Big).
\end{aligned}
\end{equation}
From \eqref{Pijk1}--\eqref{canshow} and \eqref{Em-1}, it follows that the exact formula for $P_n^{(k)}(\sigma^{-1}_i<\sigma^{-1}_j)$ is
\begin{equation}\label{thm2result2full}
\begin{aligned}
&P_n^{(k)}(\sigma^{-1}_i<\sigma^{-1}_j)=\\
&\frac12+\frac1{2n(n-1)(n-2)}\Big(\left(2(k-1)(j-i)+k^2-3k+1\right)n-2(k^2-k-1)(j-i)\Big)+\\
&\frac{\frac{(-1)^{n-k}}{(n-k)!}}{\sum_{l=0}^{n-k}\frac{-1)^l}{l!}}\thinspace\frac1{2n(n-1)(n-2)}\times\\
&\Big(kn^3-(5k+1)n^2+(-k^3(3-2i)k^2+9k+2(j-i)+1)n+\\
&2ik^3+2(i-2)k^2-2(j-i+2)k-2(j-i)\Big).
\end{aligned}
\end{equation}

\end{document}